\documentclass[11pt]{article}
\usepackage{amsmath}
\usepackage{amsthm}
\usepackage{amssymb}

\newtheorem{theo}{Theorem}[section]

\newtheorem{prop}[theo]{Proposition}

\newtheorem{conj}[theo]{Conjecture}

\newcommand{\CC}{{\cal C}}

\newcommand{\eps}{{\varepsilon}}
\begin{document}
\date{}

\title{
Hitting all maximum independent sets
}

\author{Noga Alon
\thanks
{Princeton University,
Princeton, NJ 08544, USA and
Tel Aviv University, Tel Aviv 69978,
Israel.
Email: {\tt nalon@math.princeton.edu}.  
Research supported in part by
NSF grant DMS-1855464, BSF grant 2018267
and the Simons Foundation.}
}

\maketitle
\begin{abstract}
We describe an infinite family of graphs $G_n$, where
$G_n$ has $n$ vertices,
independence number at least $n/4$, and no set of less than
$\sqrt{n}/2$ vertices intersects all its maximum independent sets.
This is motivated by a question of Bollob\'as, Erd\H{o}s and Tuza,
and disproves a recent conjecture of Friedgut, Kalai and Kindler.
Motivated by a related  question of the last authors, we show
that for every graph $G$ on $n$ vertices with independence
number $(1/4+\eps)n$, the average independence number of an induced
subgraph of $G$ on a uniform random subset of the vertices is at
most $(1/4+\eps-\Omega(\eps^2)) n$. 
\end{abstract}

\section{Background and results}
The following conjecture appears in a recent paper of Friedgut,
Kalai and Kindler.
\begin{conj}[\cite{FKK}, Conjecture 3.1]
\label{c11}
For every $\alpha \in (0,1/2)$ there exists $k$ and $\tau>0$ 
such that if $G$ is a graph on $n$ vertices with maximum
independent set of size $\alpha n$, then there exist pairwise
disjoint subsets of vertices $A_1,A_2, \ldots ,A_r$ in $G$
such that 
\begin{enumerate}
\item
$|A_i|=k $ for all $i$.
\item
$|\cup_{i=1}^r A_i| \geq \tau n$
\item
Every maximum independent set in $G$ intersects every set $A_i$.
\end{enumerate}
\end{conj}
Here we show that this is false for every fixed $\alpha \in
(0,1/2)$ even if the requirement (2) is omitted, but that the
assertion does hold for any $\alpha >1/2$. We also discuss several
related problems.

For a graph $G=(V,E)$ let $h(G)$ denote the minimum cardinality of
a set of vertices that intersects every maximum independent set of
$G$. Bollob\'as, Erd\H{o}s and Tuza (see \cite{Er}, page 224, or
\cite{CG}, page 52) raised the following conjecture.
\begin{conj}[\cite{Er}, \cite{CG}]
\label{cbet}
For any positive $\alpha$,
if the size $\alpha(G)$ of a
maximum independent set in an $n$-vertex graph $G$ 
is at least $\alpha n$, then $h(G)=o(n)$.
\end{conj}
They formulated the conjectures for cliques and not for
independent sets; replacing $G$ by its complement leads to the
formulation above.

\begin{theo}
\label{t12}
For every positive integer $k$ there is a graph $G=G_k$ with
$n=2k(2k-1)$ vertices, independence number $\alpha(G)=k^2 (> n/4)$,
and $h(G)=k+1 (> \sqrt {n}/2).$
\end{theo}
\begin{theo}
\label{t13}
For any positive integers $m$ and $t$, where $m$ is even and $4t^2
\leq m$,
there is a graph $G=G_{m,t}$ on $n=2^m$ vertices with independence number 
$\alpha(G)=\sum_{i=0}^{m/2-t} { m \choose i}$ and
$h(G)=\Theta(t^2)$.
\end{theo}
The above two results provide counterexamples to
Conjecture \ref{c11}. On the other hand we observe that an old
result of Hajnal implies that the
assertion of the conjecture (with $k=1$ and $\tau=2\alpha-1$) 
does hold for any $\alpha >1/2$. Theorem \ref{t13} also settles 
the final open problem raised by Dong and Wu in \cite{DW}.

The graphs establishing the assertion of Theorem \ref{t12} are
regular. It turns out that for regular graphs (of any degree)
the estimates in this theorem are nearly tight, as stated 
in the next proposition.
\begin{prop}
\label{p14}
For any fixed $\eps>0$ and any regular graph $G$ with $n>n_0(\eps)$ 
vertices
satisfying $\alpha(G)  \geq (1/4+\eps)n$, the parameter $h(G)$
satisfies $h(G) < (1/\eps) \sqrt {n \log n}+1$.
\end{prop}

Conjecture \ref{c11} is motivated by another conjecture suggested
in \cite{FKK}.
For a graph $G$ on $n$ vertices 
and independence number $\alpha(G)=\alpha n$,
let $\alpha'n=\alpha'(G)n$ 
denote the average value of the independence number of
the induced subgraph of $G$ on a uniform random set of
vertices. 
\begin{conj}[\cite{FKK}, Conjecture 2.9]
\label{c15}
For any $\alpha \in (0,1/2)$ there is an $\eps=\eps(\alpha)>0$
so that for every graph $G$ with $n$ vertices and 
independence number $\alpha n$,
$\alpha'(G) \leq \alpha-\eps(\alpha)$.
\end{conj}
The following result shows that the above is true for any
$\alpha > 1/4$. The case of positive $\alpha$ close to $0$ appears
to be significantly more difficult (and interesting).
\begin{theo}
\label{t16}
Let $G=(V,E)$ be a graph with $n$ vertices and independence
number $\alpha(G) = (1/4+\eps)n$, where $\eps>0$ satisfies
$\eps<1/4$. Then 
$\alpha'(G) \leq 1/4+\eps-\eps^2/3$.
\end{theo}

The constant $1/3$ above can be easily improved and we make no attempt to
optimize it here.

For regular graphs we can show that the assertion of the conjecture
holds for all fixed $\alpha>1/8$.
\begin{prop}
\label{p17}
For any $\eps>0$  there is some $g(\eps)>0$ so that the following
holds.
Let $G=(V,E)$ be a regular graph with $n$ vertices and independence
number $\alpha(G) = (1/8+\eps)n$.
Then 
$\alpha'(G) \leq 1/8+\eps-g(\eps)$.
\end{prop}

The proofs appear in the next section. All logarithms throughout
the paper are in
base $2$, unless otherwise specified. To simplify the presentation
we omit all floor and ceiling signs whenever these are not crucial.

\section{Proofs}

\subsection{Constructions}

\noindent
{\bf Proof of Theorem \ref{t12}:}\, 
The graph $G=G_k$ is the shift graph described as follows.
Put $K=\{1,2,.., 2k\}$. The  set of 
vertices of $G_k$ is the set of all
ordered pairs $(i,j)$ with $i \neq j$ and $i,j \in K$. Thus the 
number of vertices is 
$n=2k(2k-1)$. Two vertices $(a,b)$ and $(c,d)$ are adjacent
if $b=c$ or $d=a$.
Note that the vertices can be viewed as all
directed
edges of the complete directed graph on $K$, where 
two are adjacent iff they form
a directed path of length $2$. It is easy to check 
that the maximum independent sets of this graph are of size
$k^2$. Indeed, for every partition of $K$ into two disjoint 
parts $S$ and $T$
of equal cardinality, the set of all
pairs $(s,t)$ with $s \in S, t \in T$ is a maximum independent
set, and these are all the maximum independent sets.
Any set $H$ of at most $k$ vertices of $G$ can be viewed as
$k$ directed edges of the complete graph on $K$. Let $S$
be a set of $k$ points in $K$ that does not contain the 
head of any of these $k$ directed edges, and put $T=K-S$.
Then the maximum independent set consisting of all
pairs $(s,t)$ with $s \in S, t \in T$ does not intersect $H$.
Therefore $h(G) \geq k+1$. This is tight as shown by 
a set of pairs forming a directed cycle of length $k+1$ in the
complete directed graph on $K$. \hfill $\Box$
\vspace{0.2cm}

\noindent
{\bf Proof of Theorem \ref{t13}:}\, 
Let $G=G_{m,t}$ be the graph whose vertices are all binary vectors
of length $m$, where two are adjacent iff the Hamming distance
between them exceeds $m-2t$. Note that this is the Cayley graph
of $Z_2^m$ with respect to the set of all vectors of Hamming weight
at least $m-2t+1$. This graph contains as an induced subgraph
the Kneser graph $K(m,m/2 +1 -t)$.
By an old result of Kleitman \cite{Kl},
the independence number of this graph is
exactly $\sum_{i=0}^{m/2-t} {m \choose i}$.
The maximum independent sets are the $2^m$ Hamming
balls of radius $m/2-t$ centered at the vertices of $G$. Any set of
vertices that hits all these independent sets forms a covering code
of radius $m/2-t$ in $Z_2^m$. By using known results about covering
codes  in this range of the parameters 
it is not
difficult to prove that the minimum possible size of such a set is
$\Omega(t^2)$. Indeed, viewing the vectors of the covering code
as vectors with $\{-1,1\}$ coordinates, if their number is $T$ then 
by a known result in Discrepancy Theory (see, e.g., \cite{AS}, 
Corollary 13.3.4), there is a $\{-1,1\}$ vector whose inner product 
with all members of the code is in absolute value at most
$12 \sqrt T$. If $12 \sqrt T < 2t$ this gives a vector 
whose Hamming distance from any codeword is larger than $m/2-t$,
contradicting the assumption. This shows that the size of the code
is at least $\Omega(t^2)$. 
This is tight up to the hidden constant in the
$\Omega$-notation as can be shown by a random construction of
vectors of length $\Theta(t^2)$, extending each such vector in two
complementary ways on the remaining coordinates, or by taking the
rows of a Hadamard matrix of order $\Theta(t^2)$ and their
inverses, extending them in the same way.
Note that the fact that the Kneser graph
$K(m,m/2+1 -t)$ is a subgraph of $G$ also
implies a lower bound of $2t$ for the size of the hitting set
(as the Hamming balls  of radius $m/2-t$ centered in the points of
the hitting set cover all points, providing a proper coloring of the
Kneser graph), but the bound obtained this way is weaker than the
tight $\Theta(t^2)$ bound. \hfill  $\Box$

\subsection{Induced subgraphs on random subsets}

\noindent
An old result of Hajnal \cite{Ha} (see also \cite{Ra}) asserts
that for every graph $G$ the cardinality of the intersection of all 
maximum independent sets plus the cardinality of the union of all
these sets is at least $2 \alpha(G)$.  If $\alpha(G) = \alpha n$
where $\alpha>1/2$ and $n$ is the number of vertices of $G$,
this implies that there is a set of at least $(2\alpha
-1)n$  vertices contained in all maximum independent sets. 
(The result in \cite{Ha} is formulated in terms of cliques, and not 
in terms of independent sets, but this is clearly equivalent as
shown 
by replacing the graph with its complement). 
This shows that the assertion of Conjecture
\ref{c11} holds for $\alpha>1/2$ (with $k=1$ and $\tau=2\alpha-1$).

Using Hajnal's result we next describe the proof of Theorem
\ref{t16} obtaining a 
(modest) progress in the study of Conjecture \ref{c15}.
\vspace{0.2cm}

\noindent
{\bf Proof of Theorem \ref{t16}:}\,  
Without loss of generality we may assume that
$n$ is arbitrarily large, as we can replace $G$ by a union
of many vertex disjoint copies of itself and use linearity of
expectation. Assuming $n$ is large, almost every random subset of
vertices is of cardinality $(1/2+o(1))n$, hence it suffices to
show that for almost every set $W$ of $m=(1/2+o(1))n$ vertices,
the independence number of the induced subgraph of $G$ on $W$ is
smaller than $(1/4+\eps-\eps^2/2)n$. Construct the random set $W$
of size $m$ by  removing from $G$ vertices, one by one. Starting
with $V=V_0$, let $V_{i+1}$ be the set obtained from $V_i$ by 
removing a uniform random vertex of $V_i$. The set $W$ is thus
$V_{n-m}$. Let $G_i$ be the induced subgraph of $G$ on $V_i$.
Call a step $i$, $1 \leq i \leq
n-m$ of the random process above 
successful if either the independence number of 
$G_{i-1}$ is already smaller than $(1/4+\eps-\eps^2/2)n$
(note that in this case this will surely be the case in
the final graph $G_{n-m}$), or the independence number
of $G_{i}$ is strictly  smaller than that of $G_{i-1}$.
Put $i_0=(1/2-\eps)n$ and consider the graph 
$G_{i_0}$. For any
$i >i_0$, the number of vertices of $G_{i-1}$ is at most
$(1/2+\eps)n$. If its independence number is  smaller than
$(1/4+\eps-\eps^2/2)n$ then, by definition, step number $i$ is 
successful. Otherwise, by the result of Hajnal mentioned above,
the number of vertices of $G_{i-1}$ that lie in all the maximum
independent sets in it is at least $(\eps-\eps^2)n$. Since $\eps <1/4$
this is a fraction of at least  $\eps$ of the vertices of
$G_{i-1}$. Therefore, in this case,
the probability that the next chosen vertex lies in all maximum
independent sets of $G_{i-1}$ is at least $\eps$. 
We have thus shown that for every $i$ satisfying $i_0<i  \leq n-m$
the probability that
step number $i$ is successful is at least $\epsilon$. Therefore,
the probability that there are at least $\eps^2 n/2$ successful
steps during the $n-m-i_0=(\eps-o(1))n$ steps starting with
$G_{i_0}$ until we reach $G_{n-m}$  is at least the probability
that a binomial random variable with parameters
$(\eps-o(1))n$ and $\eps$ is at least $\eps^2 n/2$. This
probability is $1-o(1)$ for any fixed positive $\eps$ as 
$n$ tends to infinity. Since having that many successful steps 
ensures that the independence number of the induced subgraph
of $G$ on $W$ is at most $(1/4+\eps-\eps^2/2)n$, this completes
the proof.  \hfill $\Box$

\subsection{Regular graphs}
In  the proofs of Propositions \ref{p14} and \ref{p17}  we apply
the following early version of the  container theorem of \cite{BMS} and
\cite{ST}.
\begin{theo}[c.f. \cite{AS}, Theorem 1.6.1] 
\label{t21}
Let $G=(V,E)$ be a $d$-regular graph on $n$ vertices  and let
$\delta>0$ be a positive real. Then there is a collection
$\CC$ of subsets of $V$ of cardinality
$$
|\CC| \leq \sum_{i \leq n/\delta d} { n \choose i} 
$$
so that each $C \in \CC$ is of size at most
$\frac{n}{\delta d} +\frac{n}{2-\delta}$ and every independent set
in $G$ is fully contained in a member $C \in \CC$.
\end{theo}
\vspace{0.2cm}

\noindent
{\bf Proof of Proposition \ref{p14}:}\, 
Let $G=(V,E)$ be a $d$-regular graph on $n$ 
vertices with independence number at least $(\frac{1}{4}+\eps)n$,
and assume that
$n$ is sufficiently large as a function of $\eps$.
The closed neighborhood
of any vertex of $G$ intersects  every maximum independent set
of $G$, implying that $h(G) \leq d+1$.  If 
$d \leq (1/\eps) \sqrt {n \log n}$ this
implies the desired result, hence we may and will assume that
$d$ is larger. By Theorem \ref{t21} with $\delta=\epsilon$ 
there is a collection $\CC$ of at most
$$
\sum_{i \leq \sqrt n/ \sqrt {\log n} }
{ n \choose i} \leq 2^{\sqrt {n \log n} }
$$
subsets of $V$, each of size at most
$$
\frac{n}{\sqrt {n  \log n}} +\frac{n}{2-\eps}<
(\frac{1}{2}+\eps)n
$$ 
so that every independent set of $G$ is fully
contained in one of them. 

Let $X$ be a random set of $\frac{1}{\eps} \sqrt {n \log n}$
vertices chosen uniformly (with repetitions) among all
vertices of $G$. Fix a container $C \in \CC$. By Hajnal's
result there are at least $\eps n$ vertices contained in all 
maximum independent sets of $G$ that are contained in $C$. The
probability that $X$ does not contain any of these vertices
is at most 
$$(1-\eps)^{\frac{1}{\eps} \sqrt {n \log n}}
\leq e^{-\sqrt {n \log n}}.
$$
The desired result follows by applying the 
union bound over all $C \in \CC$.
\hfil $\Box$
\vspace{0.2cm}

\noindent
{\bf Proof of Proposition \ref{p17}:}\, 
Let $G=(V,E)$ be as in the proposition, where $|V|=n$.
As before we may assume without loss of generality that
$n$ is sufficiently large as a function of $\eps$.
Without trying to optimize the function $g(\eps)$, let $d$ denote
the degree of regularity of $G$. 
Note that $G$ contains  a set $S$ of  at least $n/(d^2+1)$ vertices 
no two of which are adjacent or have a common neighbor.
Let $W$ be a uniform random set of vertices of $G$. If the complement
of $W$  fully contains the closed neighborhoods of $s$ vertices
of $S$, then  the independence number of the induced subgraph of
$G$ on $W$ is at most $\alpha(G)-s$. The random variable counting
the above number $s$ is a Binomial random variable with
expectation $|S|/2^{d+1} \geq \frac{n}{(d^2+1)2^{d+1}}$. 
Thus if, say, $d$ is at most $50/\eps^4$ we get that the expected
independence number of an induced subgraph of $G$ on a uniform
random set of vertices is at most 
$$
\alpha(G) -\frac{n}{(d^2+1)2^{d+1} }
$$
supplying a lower bound (of the form $2^{-\Theta(\eps^{-4})}$) 
for $g(\eps)$. We thus may and will assume that 
$d \geq \frac{50}{\eps^4}$. By Theorem \ref{t21}
with $\delta=\eps$ there is a collection $\CC$ of subsets of 
$V$, satisfying
$$
|\CC| \leq \sum_{i \leq \eps^3n/50} {n \choose i} 
\leq 2^{H(\eps^3/50) n},
$$ 
where $H(x)=-x \log x -(1-x) \log (1-x)$ is the binary entropy
function. Each member $C$ of $\CC$ is of size at most
$$
\frac{n}{\delta d} + \frac{n}{2-\delta} 
\leq \frac{\eps^3 n}{50}+ \frac{n}{2-\eps} < (\frac{1}{2}+\eps)n
$$
and every independent set of $G$ is contained in a member $C \in
\CC$.

As in the proof of Theorem \ref{t16} we can generate a random
subset $W$ of $V$ by omitting vertices one by one, starting with
$V$. Since $n$ is large almost all sets $W$ are of size $n/2+o(n)$.
Moreover, for almost all of them the size of $W \cap C$ deviates
from $|C|/2$ by at most, say, $\frac{\eps}{100} n$
for all $C \in \CC$, provided $\eps$ is sufficiently small. 
It suffices to show that with high probability the independence
number of the induced subgraph on $W$ is
at most, say, $\alpha(G)-2g(\eps) n$.
Since every
independent set is contained in at least one of the members
$C$ of $\CC$ it suffices to show that with high probability 
the independence number of the induced subgraph of $G$ on $W \cap
C$ is at most $\alpha(G)-2g(\eps)n$ for every $C \in \CC$. 
Fix $C \in \CC$. Without loss
of generality its size is at least $n/8$ (since otherwise it
cannot contain a large independent set at all). Recall that
$|C| \leq (1/2+\eps)n$. If the independence number of the induced
subgraph of $G$ on $C$ is smaller than $(1/4+\eps)|C|$ then so is
the independence number of the induced subgraph on $W \cap C$,
and this is smaller than $\alpha(G) -0.1 \eps n$, implying the
desired result.
Otherwise, as in the proof of Theorem \ref{t16}, in the random 
process that omits vertices of $C$ one by one to get
$W \cap C$, the number of times the independence number
drops dominates stochastically a binomial random variable 
with parameters $\frac{\eps}{2} |C|$ and $\eps$. By the standard
estimates for Binomial distributions (c.f., e.g., \cite{AS},
Theorem A.1.13), the probability this  variable is less
than half its expectation is at most
$$
e^{- \eps^2 |C|/16} \leq e^{-\eps^2 n/128}.
$$
By the union bound over all $C \in \CC$ the probability this 
happens even for a single $C \in \CC$ is at most
$$
2^{H(\eps^3/50)n} \cdot e^{-\eps^2 n /128}
$$
which, for small $\eps$, tends to $0$ as $n$ tends to infinity.
This shows that in this case  ($d \geq \frac{50}{\eps^4}$ ),
with high probability the independence number
of the induced subgraph of $G$ on $W \cap C$ is smaller than
$\alpha(G)$ by at least, say, $\eps^2  n/40$, completing the 
proof of the proposition.  \hfill $\Box$

\section{Remarks}
\begin{itemize}
\item
Conjecture \ref{cbet} remains open for $n$-vertex graphs
with independence number at most $n/2$ and for such regular graphs
of independence number at most $n/4$. Similarly,
Conjecture \ref{c15} remains open for 
$n$-vertex graphs
with independence number at most $n/4$ and for such regular graphs
with independence number at most $n/8$. Both conjectures
appear to be significantly more difficult for graphs with
independence number $\beta n$ when $\beta>0$ is a fixed  small
positive real.  
\item
A conjecture I raised more than ten years ago motivated 
by some of the results in \cite{AHLSW}
is that the chromatic number of the graph $G_{m,t}$ described
in the proof of Theorem \ref{t13}, where $4t^2 \leq m$,
is $\Theta(t^2)$. This has been mentioned in several lectures,
see, for example, \cite{Al}. By the arguments described in the
proof of Theorem \ref{t13} this chromatic number is at least
$2t$ and at most $O(t^2)$.
\end{itemize}

\noindent
{\bf Acknowledgment}
I thank Ehud Friedgut for helpful discussions and Zichao Dong and
Zhuo Wu for telling me about \cite{DW}.

\end{document}